\documentclass[3p]{elsarticle}
\usepackage{amssymb,amsmath}

\usepackage{graphicx}
\usepackage{cancel}
\usepackage[ruled]{algorithm2e}

\usepackage[colorinlistoftodos,textwidth=4cm,shadow]{todonotes}
\newcounter{Ivan}

\newcounter{Igor}

\newcounter{Denis}

\date{}
\title {\textbf{What Lies Beneath the Surface: Topological-Shape Optimization With the Kernel-Independent Fast Multipole Method}}

\author{Igor Ostanin \footnote{Corresponding author, tel: +79150174677, e-mail:i.ostanin@skoltech.ru}
}
\author{Ivan Tsybulin}
\author{Mikhail Litsarev}
\author{Ivan Oseledets}
\author{Denis Zorin}

\address{Skolkovo Institute of Science and Technology, Nobel St. 3, Moscow, Russia}

\begin{document}

\begin{abstract}

The paper presents a new method for shape and topology optimization based on an efficient and scalable boundary integral formulation for elasticity. To optimize topology, our approach uses iterative extraction of isosurfaces of a topological derivative. The numerical solution of the elasticity boundary value problem at every iteration is performed with the boundary element formulation and the kernel-independent fast multipole method. Providing excellent single node performance, scalable parallelization and the best available asymptotic complexity, our method is among the fastest optimization tools available today. The performance of our approach is studied on few illustrative examples, including the optimization of engineered constructions for the minimum compliance and the optimization of the microstructure of a metamaterial for the desired macroscopic tensor of elasticity.       

\end{abstract}

\maketitle

\section{Introduction}

The idea of topological optimization (also known as layout optimization or structural optimization) has its roots in the classic century-old work by Mitchell \cite{Mitchell1904}, and is increasingly important due to advances in fabrication technologies making it possible to manufacture optimized shapes. In their present form, numerical techniques of topological optimization originate from the seminal paper by Bends{\o}e and Kikuchi \cite{Bendsoe1988}. Initially driven by the demands of automotive and aerospace industry, modern topology optimization techniques have applications in biomedical and electrical engineering design, architecture  and material science ( for an overview, see, for example \cite{Review, Topopt}).   

The problem of topology optimization can be stated as follows: for a given domain, boundary conditions, and a set of constraints, find a distribution of the material that minimizes a cost functional depending on the solution of the partial differential equation (PDE) of interest ( elasticity, electric or heat conductivity \textit{etc.}) in this domain. The most common  example of such an optimization is the minimization of compliance, \textit{i.e.} finding the distribution of elastic material that, for a given total weight and boundary conditions, minimizes the elastic strain energy. 

All common topology optimization methods used in commercial and academic software are based on the finite difference or finite element methods (FEM), the latter being the only practical option for complex domains.  While FEM is the most widely used and flexible approach to solving elasticity problems required for topology optimization, it has a number of common drawbacks, especially for large-scale problems: the need to discretize the whole optimization domain, and solve ill-conditioned systems. Beyond that, there are  downsides specific to topology optimization ---  emergence of spurious solutions (``checkerboard patterns'') that need to be eliminated with regularization, and mesh dependency of the optimal solutions and the corresponding cost functionals \cite{Meshdep}.  

Several recent papers ( \textit {e.g.} \cite{BEM2DMarczak, BEM3D}) demonstrated that the boundary element method (BEM) \cite{Cruse} could be used as a tool for topology optimization. In this work we explore this idea and describe a fast, scalable and numerically stable implementation of the algorithm of topology optimization based on the BEM. Our technique is shown to be free of the typical limitations of FEM formulations while demonstrating the single-node performance and parallel scalability comparable or better than the state of the art FEM solvers. We employ the recent implementation \cite{pvfmm} of the kernel independent fast multipole method (KIFMM) \cite{kifmm} in combination with BEM to address the problems of shape and topology optimization. The key features of our approach are: 

\begin{itemize}
 \item use of BEM formulation and discretization of the elasticity BVP; 
 \item a discrete binary method of material removal; 
 \item acceleration of the boundary element solve and topological derivative evaluation with a highly scalable fast multiple method, suporting the kernels we need for the elasticity solve and stress tensor evaluation; 
 \item extraction of the boundary as a level set of topological derivative.         
\end{itemize}

The application of the boundary integral method to elasticity and acceleration of BEM with fast multipole method are well-known, as well as the topology optimization techniques based on topological derivatives. However, no attempts were made so far to combine them into a powerful and scalable algorithm designed specifically for topological-shape optimization. Our present work addresses this problem. Below we describe the key ideas of our work and demonstrate the validity, robustness and scalability of our technique on few illustrative examples. Also, in the remainder of the paper we provide a brief overview of the related works, putting our developments into a broader context of modern trends in topological-shape optimization.

\section{Method}

We seek to solve the following problem. Consider an elastic domain $\Omega$ with the boundary $\partial \Omega $,  filled with a linear isotropic elastic material with bulk modulus $K$ and shear modulus $G$. A mixed BVP for elasticity PDEs is prescribed for this domain (here and below the indicial notation is used): 

\begin{equation} \label {elast}
\begin{split}
\sigma_{ij,j}=0,  \\
\sigma_{ij}=C_{ijkl}\epsilon_{kl}, \\
C_{ijkl} = K \cdot \delta_{ij} \delta_{kl} + G \cdot ( \delta_{ik} \delta_{jl} + \delta_{il} \delta_{jk} - 2/3 \delta_{ij} \delta_{kl}), \\
\sigma_{ij} \cdot n_j = t^N_i \mid_{\partial \Omega_N}, \\
u_i = u^D_i \mid_{\partial \Omega_D} \\
\end{split}
\end{equation}

We search for a subdomain $\omega$ that, for its given volume, minimizes the cost functional

\begin{equation} \label {functional}
\Psi=\frac{1}{2}\int_{\omega}\sigma_{ij}\epsilon_{ij}d\varOmega
\end{equation}
 
In order to address this problem, we employ a "hard-kill" approach, based on unidirectional, "hard" elimination of the material of the original domain. As a measure of sensitivity of the cost functional to material removal at a certain point of the domain, we utilize the \emph{topological derivative} (TD) \cite{Sokolmain} -- a cost of making an infinitesimal spherical cavity with a center in a given point of the domain . For the case of strain energy (compliance) cost functional and 3D linear isotropic elasticity the analytical expression for TD is available \cite{Sokol3D,Novotny3Delast}:

\begin {equation} \label{deriv}
D^{T}(p)=\frac{3}{4E}\frac{1-\nu}{7-5\nu}\left[10(1+\nu)\sigma_{ij}(p)\cdot\sigma_{ij}(p)-(1+5\nu)tr\sigma_{ij}(p)^{2}\right]
\end{equation}

here $E$ is the material's Young's modulus.

We discretize the the initial domain volume on a set of square cells. Our method works on an arbitrary connected subset of a regular grid.  Each cell of the grid is marked as \emph{filled} or \emph{empty}. The following sequence of steps is performed at every optimization iteration.

\begin{enumerate}
        \item Initialize all cells in the domain to \emph{filled}.
        \item Extract the boundary of the part of the domain filled
          with material; 
	\item Solve the boundary value problem in BEM formulation.
	\item Based on the BVP solution, compute the values of the
          topological derivative at all filled cells.
	\item Mark all cells meeting the criterion for material
          removal as \emph{empty}.
	\item Quit if desirable volume ratio of the material is
          reached, otherwise return to step 2. 	         
\end{enumerate}

Next, we discuss the the numerical solution of a BVP, criterion for
the material removal and the parallel iterative optimization procedure in
greater detail.  

\subsection{Boundary Element Formulation} 

One of the key points of our approach is the use of surface integral
equations and boundary discretization for the solution of elasticity
BVP at every optimization iteration.  We use the direct boundary integral equation (BIE) formulation for elasticity \cite{Cruse}:

\begin {equation} \label{BIE}
\frac{1}{2}u_{i}(\xi)=\int_{\varGamma}U_{ij}(\xi,x)p_{j}d\Gamma-\int_{\varGamma}P_{ij}(\xi,x)u_{j}d\Gamma
\end {equation}
 where $u_{i}(x)$ and $p_{i}(x)$  are the displacement and tractions  on the boundary of the domain, and $U_{ij}(x, \xi)$ ($P_{ij}(x, \xi)$) are corresponding fundamental solutions. For the case of a linear isotropic elastic material with the shear modulus $G$ and Poisson's ratio $\nu$ these are given by

\begin {equation} \label{disp_1}
U_{ij}(\xi,x)=\frac{1}{16\pi(1-\nu)Gr}\left((3-4\nu)\delta_{ij}+r_{,i}r_{,j}\right)
 \end {equation}
\begin {equation} \label{disp_2}
P_{ij}(\xi,x)=\frac{1}{8\pi(1-\nu)r^{2}}\left[\frac{\partial r}{\partial n}\left((1-2\nu)\delta_{ij}+3r_{,i}r_{,j}\right)-(1-2v)(r_{,i}n_{j}-r_{,j}n_{i})\right]
 \end {equation}

where $ r=\left|\xi-x\right| $. Once the solution on the boundary is found, the stress at the point inside the domain can be calculated using another integral formula: 

\begin {equation} \label{stress_1}
\sigma_{ij}(p)=-\int_{\Gamma}u_{k}(x)S_{kij}(p,x)d\varGamma+\int_{\varGamma}t_{k}(x)D_{kij}(p,x)d\varGamma
 \end {equation}
where  the fundamentall solutions $D_{kij}(p,x)$ and $S_{kij}(p,x)$ are given by 

\begin {equation} \label{stress_2}
D_{kij}(p,x)=\frac{1-2\nu}{2\pi(1-\nu)r^{2}}\left(\delta_{ki}r_{,j}+\delta_{kj}r_{,i}-\delta_{ij}r_{,k}+\frac{3}{1-2\nu}r_{,i}r_{,j}r_{,k}\right)
\end {equation}
\begin {equation} \label{stress_3}
\begin{split}
	S_{kij}(p,x)=\frac{3-6\nu}{4\pi(1-\nu)r^{3}}\left[\delta_{ij}r_{,k}+\frac{\nu}{1-2\nu}\left(\delta_{ki}r_{,j}+\delta_{kj}r_{,i}\right)-\frac{5}{1-2\nu}r_{,i}r_{,j}r_{,k}\right]\frac{\partial r}{\partial n}+ \\
	\frac{1-2\nu}{4\pi(1-\nu)r^{3}}\left[\frac{3\nu}{1-2\nu}(n_{i}r_{,j}r_{,k}+n_{j}r_{,i}r_{,k})+3n_{k}r_{,i}r_{,j}+n_{j}\delta_{ki}+n_{i}\delta_{kj}-\frac{1-4\nu}{1-2\nu}n_{k}\delta_{ij}\right]
\end{split}
\end {equation}
Below we discuss our toolkit for the fast solution of the BIE and rapid computation of the fields inside the domain.

\subsection{Numerical solution} \label{s_ns_kifmm}

The numerical treatment of the integral formulation \eqref{BIE}  requires discretization of the domain boundary only,  and therefore results in a system of linear equations that has asymptotically smaller number of  unknowns than any approach that discretizes the domain. 

However, the matrix of the resulting system is dense, and an iterative solution scheme  would require $O(N^2)$ operations per iteration, where $N$ is the number of unknowns on the boundary, which is asymptotically \emph{slower} than doing an optimal-complexity (e.g., multigrid) volume solve. 

Therefore, in order to take full advantage of the boundary integral formulation,  a fast (linear-complexity, or close) scheme for numerical solution of surface integral equations is needed.  A number of such schemes exist \cite{fmm, fmmbook, hmatrices, h2matrices, mskel, kifmm}, which make different tradeoffs between precomputation required vs. efficiency of the solve vs. generality.  For example, an $\mathcal{H}^2$-matrix method is applicable to any dense matrix, not necessarily derived from a PDE fundamental solution, but requires a relatively expensive precomputation,  while the original (analytic) FMM method requires no precomputation but a special set of translation operators needs to be derived for each kernel. We are using the \emph{kernel-independent} FMM (KIFMM) for the following reasons.  First, just as analytic FMM it requires no precomputation that depends on the surface: this is essential for our application, as the surface changes at every step.
Second, in contrast to analytic FMM, it can handle all four kernels that we need (two for the boundary integral solve, and two for the topological derivative evaluation) in an automated way: only a kernel evaluator needs to be provided. 
We use a state-of-the-art scalable implementation \cite{kifmm,pvfmm}.    

\paragraph*{Kernel-independent fast multipole method}
For completeness, we provide a brief overview of KIFMM. In order to reach linear complexity of evaluation of integrals over a surface at a large number of points simultaneously, fast multipole methods performs the following steps\cite{fmm}:

\begin{itemize}
	\item  Generation of an octree partitioning of the domain into boxes.
	\item  Fine-to-coarse tree traversal to compute compact representations of the far-field potential of a box
          (\emph{multipole expansions} for analytic FMM); these are computed hierarchically, by combining
          expansions of descendant boxes to the expansion for the parent box, using linear M2M translation
          operators. 
          Multipole expansions are used to approximate the values of the integral over all points contained in a box, 
          with the evaluation point far enough away from that box;          
	\item  a coarse-to-fine pass, that computes \emph{local expansions} for each box, that approximate the values of the integral over all points far enough away from the box inside the box. These are obtained at descendant boxes by combining the parent's local expansion (using an L2L operator) with multipole expansions of boxes that are not in the far zone of the parent, but are in the far zones of the descendants.  Multipole expansions are converted to local using M2L operators.   
	\item  At the finest level of the tree, the complete integrals are computed by adding the contributions of points in the near zone using direct summation.           
\end{itemize}

The distinguishing feature of KIFMM, compared to the original FMM method is that it does not require analytical multipole and local series expansions of underlying kernels, and analytically derived M2M, M2L, L2L operators for each kernel.
Instead, it represents the far-field (multipole) and local approximations of the integrals with a density $\phi$  defined at samples $x_i$ of an equivalent surface, so that the approximation at a point $y$  has the form $\sum_i \phi_i K(x_i,y)$, where $K(x,y)$ is the kernel of interest.  

The M2M, L2L and M2L operators needed in the algorithm, in the case of KIFMM are represented by matrices mapping density values on different equivalent surfaces to each other, and are computed automatically for each needed kernel. 

Just like the original FMM, kernel-independent FMM performs the summation of the field of $N_s$ sources at $N_t$ targets with $O(N_s + N_t)$ operations.
In our work we use a recent parallel implementation of KIFMM -- PVFMM\cite{kifmm,pvfmm}, which implements rapid evaluation of sums of the following 
\begin{equation} \label{fmm_summation}
t_{i}( \mathbf{x}_{i})=K_{ij}(\mathbf{x}_{i}, \mathbf{y}_{j},\mathbf{n}_{j})s_{j}(\mathbf{y}_{j})
\end{equation}     
$t_{i}$ is the vector of target values being computed (values of an integral at points of interest $\mathbf{x}_{i}$);  $s_{j}$ is the vector of known source  values at points $\mathbf{y}_{j}$ (solution values on the surface).  Kernel function $K_{ij}(\mathbf{x}_{i}, \mathbf{y}_{j},\mathbf{n}_{j})$  depends on both source and target coordinates, and, for double-layer kernels,  on the normals $\mathbf{n}_{j}$ that is specified at source point.
We use KIFMM for fast summation of the matrix components of kernels \eqref{disp_1}, \eqref{disp_2}, \eqref{stress_2}, \eqref{stress_3}.  PVFMM is highly optimized and extremely scalable implementation of KIFMM, it supports both intranode OpenMP standard parallelization and internode MPI standard parallelization, demonstrating excellent scalability for up to tens of thousands cores \cite{pvfmm}. 

\paragraph*{Discretization of BIE} Using collocation method and piecewise-constant approximation of tractions and displacements on the boundary \cite{Cruse}, one can discretize the equation \eqref{BIE} into the following system of linear equations:
\begin{equation} \label {num}
\left(\frac{1}{2}I+P^\varDelta \right)\cdot u=  U^\varDelta \cdot p,
\end{equation}
where 

\begin{equation} \label {coef}
\begin{split}
P_{ij}^\varDelta =\int_{S_{k}}P_{mn}(\xi_{l},x)dS_{k}, \\
U_{ij}^\varDelta =\int_{S_{k}}U_{mn}(\xi_{l},x)dS_{k}, \\
i = 3k+m, j = 3l+n, 
\end{split}
\end{equation}
where $\xi_{l}$ is \textit{l-th} collocation point and $dS_{k}$ is the area element of \textit{k-th} boundary triangle. After rearrangement of columns of matrices in \eqref{num}, we obtain the following system of linear equations, where all unknowns appear in vector $x$, while the tractions or displacements  known from boundary conditions appear in vector $y$. 

\begin{equation} \label{num2}
A x=B y.
\end{equation}     

The system matrix $A$ is neither symmetric nor positive definite. Its condition number depends on the boundary conditions and the surface geometry. We note that the coefficients of $A$ require computation of integrals in \eqref{coef},
which are singular for diagonal terms. We evaluate the non-singular off-diagonal integrals using Gaussian  quadrature for each triangle. The singular integrals over triangles are evaluated analytically.  Below, we discuss how a black-box FMM code can be used to evaluate matrix-vector product needed for solving the system \eqref{num2}.

We use  parallel implementation of GMRES algorithm \cite{Saad_GMRES} available in PETSC library \cite{petsc-user-ref,petsc-efficient} to solve this system of linear equations. Fast evaluation of matrix-vector products in \eqref{num2} is done using KIFMM, without explicit representations of matrices $A$ and $B$. 

\paragraph*{Matrix-vector products}
If the entries of the matrix are approximated using a numerical quadrature, the matrix-vector product is reduced to a sum of fundamental solutions centered at quadrature points, multiplied by displacement/traction values and quadrature weights; this is exactly the type of sums an FMM code, PVFMM in particular, is designed to compute.  However, we need to use analytic kernel integration for triangles with singularities.  To avoid adding problem-specific complexity to them FMM code, we opt for a two-pass solution. 

First, we perform rapid summation $A \cdot x$ using PVFMM. The summation \eqref{fmm_summation}  is performed with
quadrature points on triangular boundary elements $\mathbf{y}_{j}$ as source points, and values $s_{j}$ defined as the constant approximation of the solution on the triangular element, scaled by triangle's area and quadrature point's weight, and  collocation points at triangle centers used as target points. Each triangular element contains 16 quadrature points with the corresponding weights (the element is subdivided into four equal triangular parts, and the 4-node Gauss quadrature is used for each part).

Then we perform the second, \emph{local}, pass: for each target points, we subtract the inaccurate contributions from the sources corresponding to quadrature points on the triangle that contains the target point, and replace those with the analytic expressions for the singular integral over this triangle.  
 
This scheme is easily parallelizable, since all information needed at the local pass is local to a triangle. It imposes a limitation in terms of achievable model sizes - for large enough models the numerical summation of near-singular integrals over triangles neighboring to the triangle containing the collocation point becomes inaccurate; to improve accuracy, an upsampled quadrature for  triangles close to the triangle with singularity can be used. We note that this becomes an issue only when the problem size reaches tens of millions of degrees of freedom. 

\paragraph*{Computing topological derivatives}
Since the surface solution is found, the solution in stresses at internal points is found via fast summation of the kernels \eqref{stress_2}, \eqref{stress_3}. 
It worth noting, that the sampling of the points inside the domain should not necessary be uniform -- one can  use adaptive strategies of sampling points inside the domain, which reduces computational complexity of the domain computation to $O(N_{s})$, where $N_{s}$ is the number of surface targets (see the discussion in section \ref{disc}).     
      
\subsection{Parallel optimization procedure} \label{par_opt}
The optimization procedure performs the following cycle (see Algorithm \ref{Alg1}). 
We start with three-dimensional array of voxels $M(i,j,k)$ ($M(i,j,k) = 1$ corresponds to material, $M(i,j,k) = 0$ - void), at the first iteration all the voxels are material. 
We contour all the material voxels with boundary elements, and set up the surface description of BVP $\Gamma(n)$ --- a set of arrays containing  coordinates of triangle vertices, collocation points and corresponding boundary conditions, as well as the volume mesh coordinates $\Omega(i,j,k)$.
Since this part is not computationally intensive, relative to solving the elasticity problem and computing topological derivatives,  it is done serially. 

Then the surface and volume arrays $\Gamma(n)$ and $\Omega(i,j,k)$ that were generated on a master process are scattered over all MPI processes, and parallel solution of $\Gamma(n)$, as well as computation of topological derivatives inside the domain $D^{T}(i,j,k)$ and energy densities $E(i,j,k)$ is performed.
After that the field of topological derivatives $D^{T}(i,j,k)$ is gathered to a master process, and isosurface of a topological derivative is extracted by thresholding: $D_{C}^{T}=D_{min}^{T}+C(D_{max}^{T}-D_{min}^{T})$.
 All the voxels that are below the threshold are assigned to be void. 
 The parameter $C$ determines the amount of material removed at every iteration, and is chosen empirically to provide desirable rate of material removal per iteration. Typically, $C = 0.001\div 0.1$. This parameter also defines the level of details obtained in topology optimization process, therefore  $C$ can be considered as an implicit problem regularization parameter.
 After post-processing procedure excluding isolated voxels and surface irregularities we compute the value of the cost functional $\Psi$ and the ratio $\alpha$ between the current and initial number of material voxels. 
We repeat the iteration  until the target ratio $\alpha_{c}$ was not reached. 

\begin{algorithm} \label{Alg1}
	\DontPrintSemicolon 
	\KwIn{Boundary conditions, $M(i,j,k)$, $\Psi$, $C$, $\alpha_c$}
	\KwOut{$M(i,j,k)$, $\Psi$, $\alpha$}
	$\alpha \gets 1$ \\
	
	\For{Max number of cycles} {
		From $M(i,j,k)$ construct $\Gamma(n)$, $\Omega(i,j,k)$ \\
		Scatter $\Gamma(n)$,  $\Omega(i,j,k)$\\
		Parallel solution: BVP  $\Gamma(n)$, TDs $D^{T}(i,j,k)$, EDs $E(i,j,k)$  \\
		Gather  $D^{T}(i,j,k)$, $E(i,j,k)$  \\
		New $M(i,j,k)$ by thresholding $D^{T}(i,j,k)$: \\
		\For {all i, j, k}{ 
			\If{ $D^{T}(i,j,k) < D_{min}^{T}+C(D_{max}^{T}-D_{min}^{T})$ } {
				$M(i,j,k) \gets 0$ 
			}
		}		

		Post-processing of $M(i,j,k)$ \\ 
		Compute $\alpha = \frac{\sum_{ijk}M(i,j,k)}{\sum_{ijk}1}$ and $\Psi = \sum_{ijk} E(i,j,k)$ \\
		
		\If{ $\alpha<\alpha_c$ } {
			break\\
		}
	}
	\Return{ $M(i,j,k)$, $\Psi$, $\alpha$}\;
	\caption{{\sc Optimization} Performs parallel topology optimization on a uniform array of voxels}
	\label{algo:max}
\end{algorithm}

As demonstrated in the section \ref{examples}, this procedure yields reliable results. Serial operations take less than $1\%$ of total iteration time even for our largest examples, and the code in its present form was shown to be efficiently parallelizable for up to 128 cores (see the discussion below). 
For some important specific cases the general algorithm \ref{Alg1} can be substantially improved. Several  possible improvements are discussed in section \ref{disc}.

\section{Performance and scalability}

The central feature of our technique - state of the art single-node computational performance and parallel scalability. In this section we discuss the performance of our method. The computations presented in this paper were carried out on a cluster machine with nodes nx360 M4, equipped with 64 Gb RAM and Intel Xeon E5-2650 CPUs and linked with  Mellanox ConnectX-3 infiniband. Up to 8 nodes were employed in our simulations. 

In order to check the performance of a single iteration of the optimization, we consider the following benchmark problem. Consider a unit cube of the material with the elastic moduli $E = 1, \nu = 0.3$, subjected to a uniform tension (Fig. 1 (A)). The cube is discretized into volumetric cells. The refinement of discretization is defined by the number of the cells along the side of the cube $N$. Each side of the cell that belongs to a cube boundary is discretized onto 4 triangles with piecewise constant approximation of the solution. The whole boundary of the cube is therefore represented by $24 \cdot N^2$ triangles (and collocation points), $3 \cdot 16 \cdot 24 \cdot N^2$ source surface DOFs, $3 \cdot 24 \cdot N^2$ target surface DOFs. Stresses and topological derivatives are computed on the dense mesh with $3 \cdot N^3$ volume target DOFs. We perform the whole cycle of computations present in algorithm \ref{Alg1} - iterative solution of the surface BVP using GMRES algorithm, computation of stresses and topological derivatives inside the domain. GMRES convergence tolerance is set to $10^{-4}$, the algorithm converges in $7-12$ iterations, depending on the size of the model.  
Fig. 1(B,C) summarize the observed single-node performance. It is seen that our method displays linear dependence of the iteration time with respect to the number of degrees of freedom. For a wide range of sizes of the model the time required for the surface and volume solution is approximately the same. Single FMM pass (without tree construction) achieves the performance of $30k$ DOF/s on a single core for the kernels \eqref{disp_1}, \eqref{disp_2}, and $2k$ DOF/s on a single core for the kernels \eqref{stress_1}, \eqref{stress_2}. The former is significantly higher than the performance of the state of the art FEM solvers, while the latter is comparable \cite{petsc_fastest}.

\begin{figure}
	\includegraphics[width=16cm]{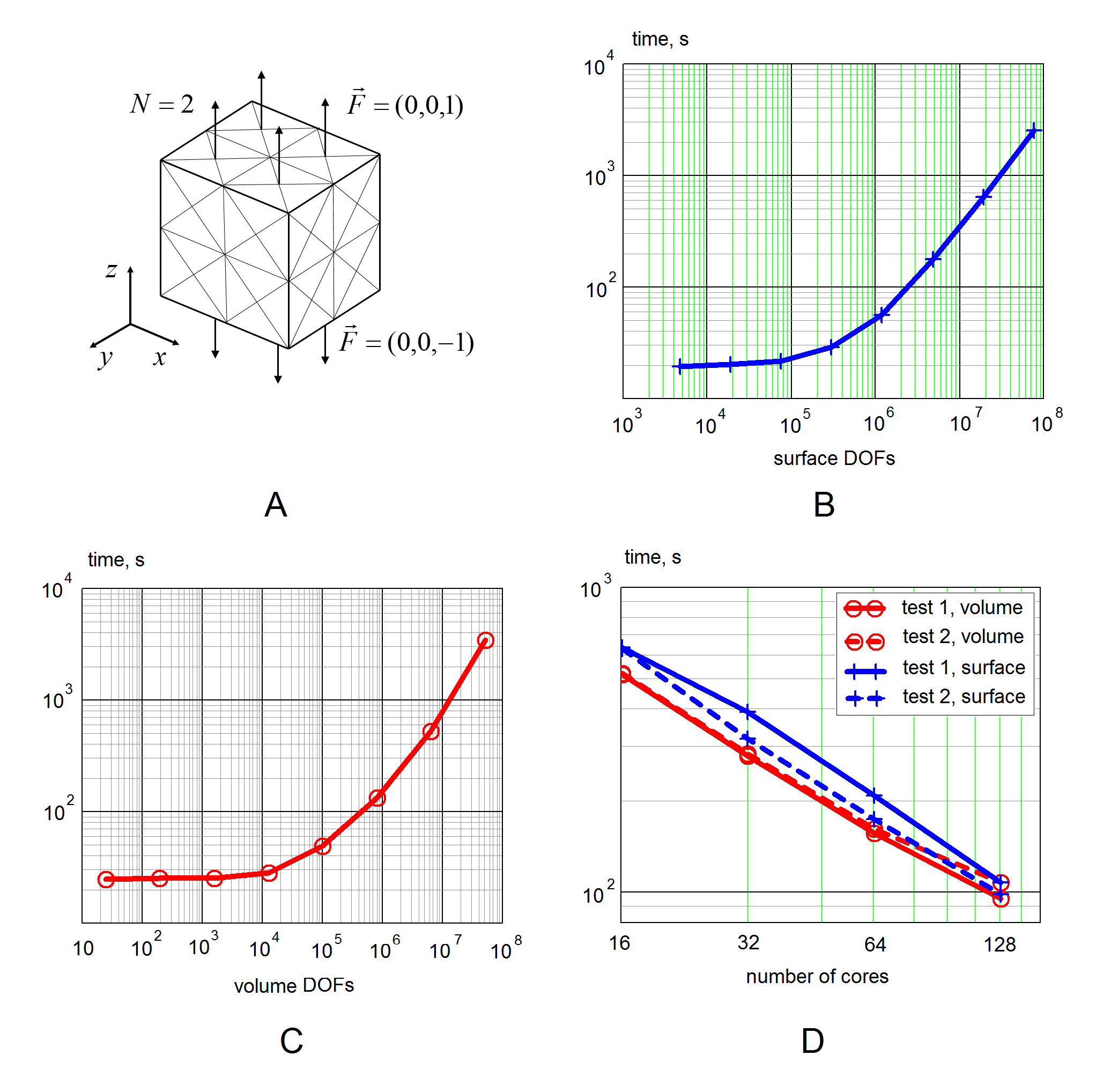}
	\protect\caption{ (A) The benchmark model ($N = 2$). (B,C) Single node performance test: (B) time to solve the surface BVP vs. the the number of the surface degrees of freedom, (C) time to compute the volume solution vs the number of the volume degrees of freedom. (D) The time required to solve the test problem with $N = 128$ on $16$, $32$, $64$ and $128$ cores of a cluster machine for two different OMP/MPI patterns. }
\end{figure}

Fig. 1(D) summarizes the strong parallel scalability. Our hybrid code supports both OMP and MPI parallelization. We performed parallel solution of the problem with $N = 128$ using $16$, $32$, $64$ and $128$ cores of the cluster machine. Two tests were performed with different OMP/MPI patterns: test 1 was carried out with the number of MPI process corresponding to a number of nodes, and 16 OMP threads running at every machine. Test 2 was carried out with 4 MPI processes per node, each running 4 OMP threads. As can be seen, the single iteration performance scales well, demonstrating approximately 6-fold performance increase in 8-node simulation, as compared to a single-node performance. 

In the following section we consider parallel optimization of benchmark optimization problems.

\section{Numerical examples} \label{examples}

This section presents few benchmark examples of topological-shape optimization with our technique. These include compliance minimization of elastic structures, as well as the optimization of a periodic cell of a metamaterial for maximum bulk modulus.
 
\begin{figure}
	\includegraphics[width=16cm]{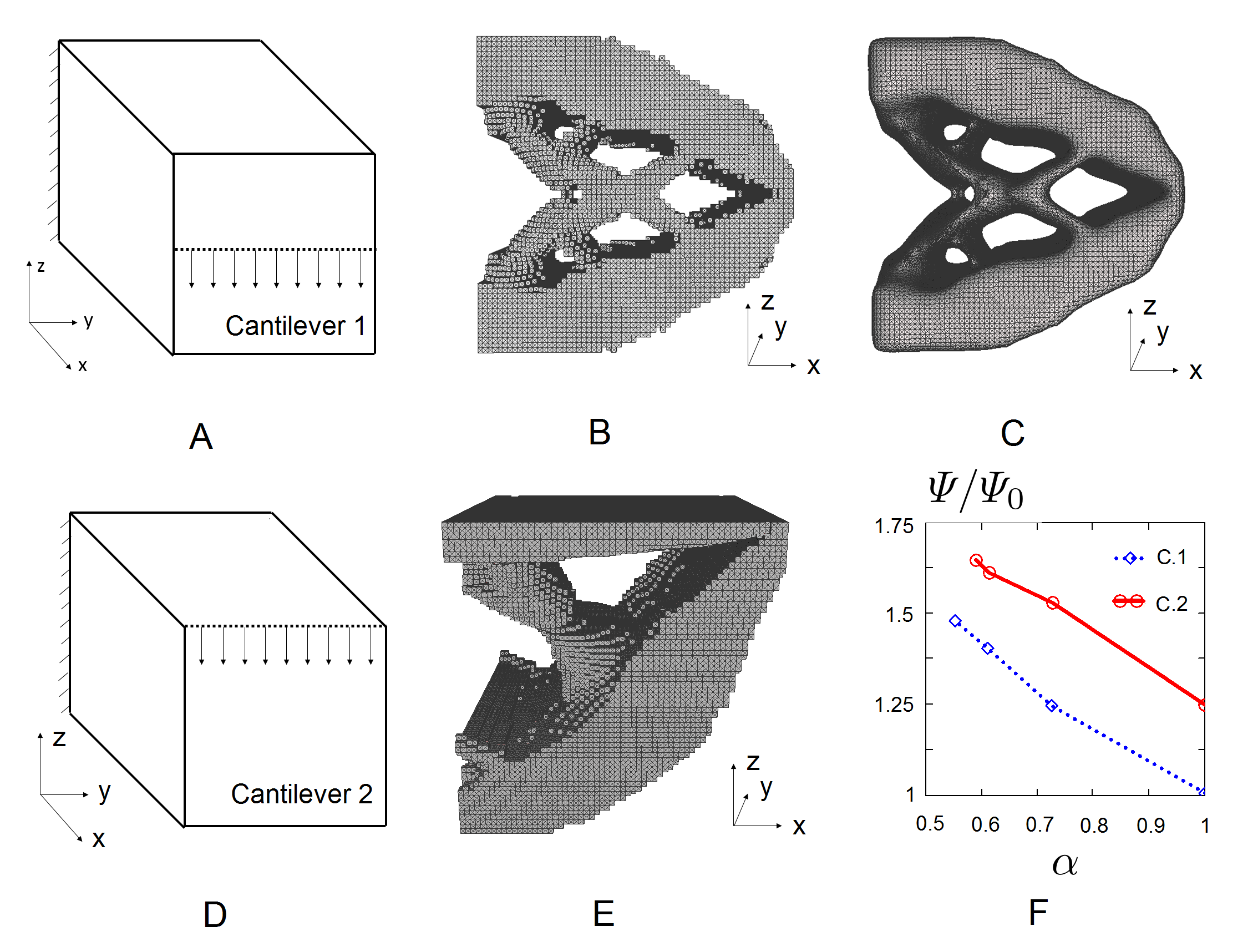}
	\protect\caption{First configuration of a cantilever support: (A) Initial BVP, (B) optimal solution for the volume fraction $\alpha_{c} = 0.6$, (C) surface after Laplacian smoothing. Second configuration of a cantilever support: (D) Initial BVP, (E) optimal solution for the volume fraction $\alpha_{c} = 0.6$.  (F) Normalized cost functional $\Psi/\Psi_0$ as a function of the current volume ratio $\alpha$ (Configurations 1 and 2).}
\end{figure}

\subsection{Cantilever supports}

We test our method on a standard example used by many authors to validate  compliance energy optimization methods. We start with a unit cube with one of its sides fixed and two different load configurations applied to the side opposite to a fixed one (Fig.1(A, D)).  The material properties are $G=1$, $\nu = 0.3$.  The initial model contains $64^3$ voxels, $0.79 \cdot 10^6$ volume DOFs and $4.72 \cdot 10^6$ surface DOFs. We use $\alpha_c = 0.6$. Computations were performed  on a cluster with up to 128 cores used. For both loading scennarios, each optimization iteration took about 7 minutes. This time was mostly determined by  relatively slow GMRES convergence, due to singularity of the solutions of BVPs with mixed boundary conditions: $100-200$ iterations were required.  
Fig. 1(B,E) show the solutions obtained after three iterations. The level of details in the final solution depends on the threshold of the topological derivative and the number of iterations. However, both obtained solutions are in  agreement with 2D and 3D solutions of similar problems obtained earlier \cite{Ostanin, Ostanin_wit, Topopt}. The quality of the surface of the optimal configuration is improved with Laplacian smoothing \cite{Laplassian} post-processing step (Fig. 1(C)). Fig. 1(F) gives the evolution of the cost functional  $\Psi$, normalized by the initial value of the cost functional for an intact cube $\Psi_0$, as a function of the current volume ratio for both examples. We can see that the symmetric configuration appears significantly stiffer than the non-symmetric one with the same material volume fraction. 
Fig. 2(A, B) demonstrate the evolution of the field of topological derivatives in a symmetry plane cross section during optimization process.

\begin{figure}
\includegraphics[width=16cm]{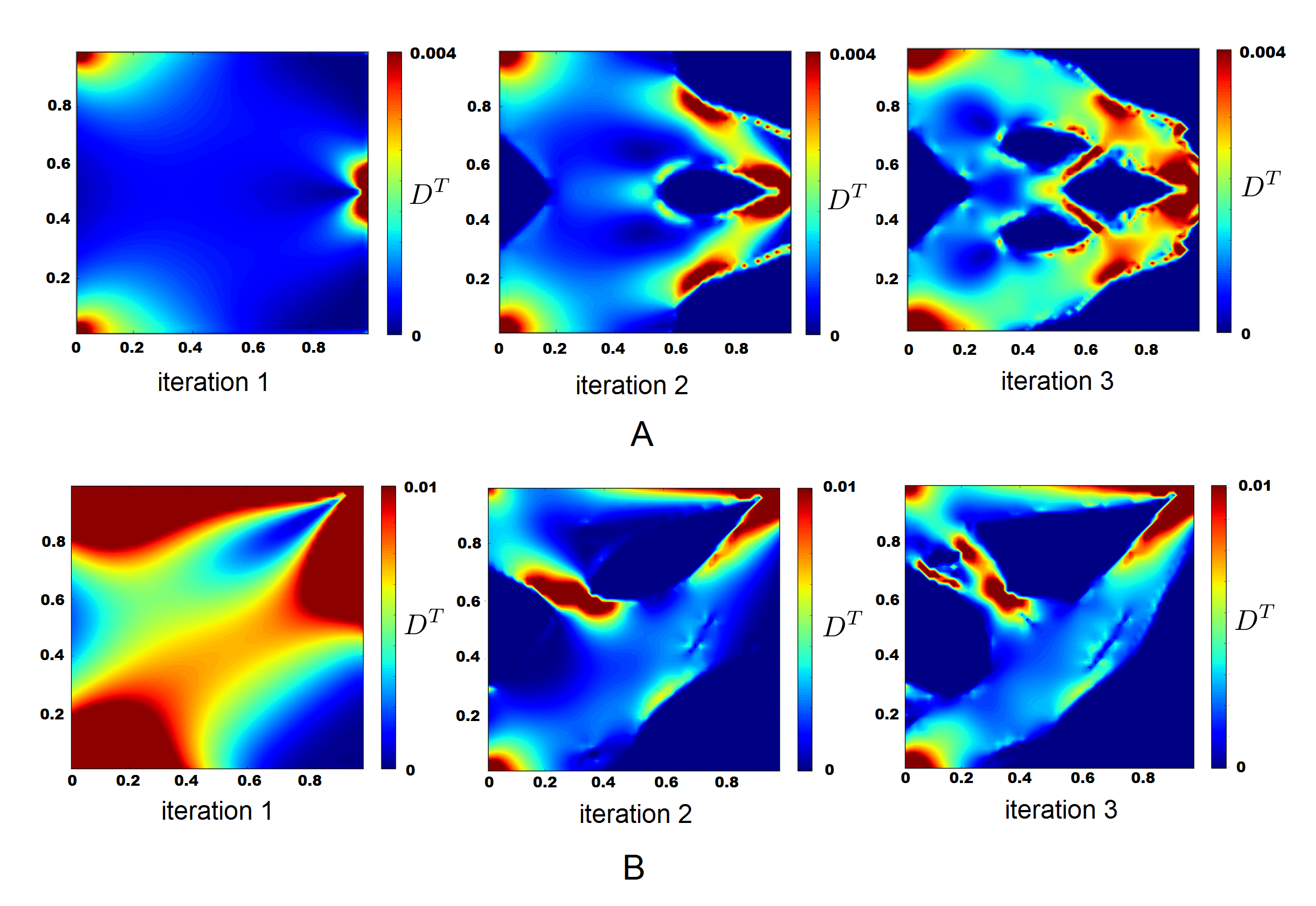}
\protect\caption{Field of topological derivatives at the $xz$ cross section passing through symmetry plane, for three optimization iterations: (A) Example 1, and (B) Example 2.}
\end{figure}

\subsection{Truss under torsion}

The following example demonstrates the compliance minimization problem in pure Neumann formulation. The initial unit cube volume of the material with $G=1, \nu = 0.3$ is subjected to a torsional loading, imposed as eight concentrated forces applied at cube vertices: $\vec{F}_{1} = (0,-1,1)$, $\vec{F}_{2} = (0,1,1)$, $\vec{F}_{3} = (0,1,-1)$,  $\vec{F}_{4} = (0,-1,-1)$,   $\vec{F}_{5} = (0,1,-1)$, $\vec{F}_{6} = (0,-1,-1)$, $\vec{F}_{7} = (0,-1,1)$, $\vec{F}_{8} = (0,1,1)$. We optimize the shape for the minimum compliance.  The model had $100^3$ voxels, each iteration took about 20 minutes on $32$ cores (again, this time was mostly conditioned by slow GMRES convergence - solution of each BVP took about 50-70 GMRES iterations). Figures 4(B) and (C) give the surfaces obtained in first and tenth iterations of the optimization process. We can see that the distinctive pattern with X-shaped frames has emerged after the first iteration and the subsequent iterations resulted in only small incremental changes in shape. 
Figure 4(D) provides the value of the normalized cost functional $\Psi/\Psi_0$ as a function the material volume fraction $\alpha$. As we can see, for the chosen value of $C$ leads to a fast convergence (1-2 iterations) in terms of the functional value.
This situation is typical for the optimization process with relatively large threshold value $C$. As we will demonstrate in the next example, the first iteration of topology optimization can immediately lead to a good optimization results and the functional values.

\begin{figure}
	\includegraphics[width=16cm]{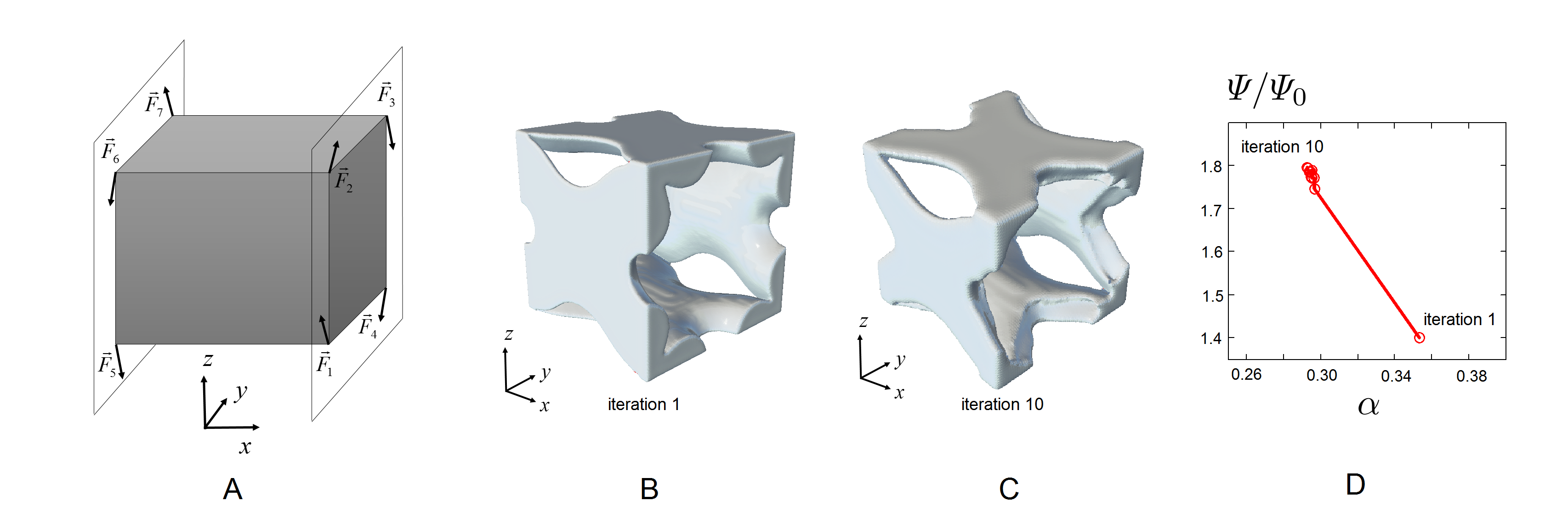}
	\protect\caption{ Truss under torsion. (A) Problem geometry and loading. Shape of the truss after the first (B) and tenth (C) iteration of the optimization process. (D) The dependence of the normalized cost functional $\Psi/ \Psi_0$ on the material volume fraction $\alpha$ .}
\end{figure}

\subsection{Periodic cell of a metamaterial} \label{periodic}

This example demonstrates how our technique can be applied to the design of a periodic cell of an elastic metamaterial. The theory of application of topological derivatives to the optimization of a periodic cell of a metamaterial is well-developed \cite{NovotnySokolovskyBook,periodicCell1,periodicCell2}.  We consider the example of the design of the metamaterial cell maximizing bulk modulus (the general case with arbitrary target homogenized properties is analogous and will be considered in a separate publication.) Our formulation of the  elasticty BVP for a periodic cell uses pure Neumann formulation, realization of integration over periodic domains in FMM, and the principle of superposition. We start with an external problem for a traction-free cavity in an infinite elastic medium subject to a homogeneous state of stress. Total stress field can be represented as a superposition of two elastic fields: a homogeneous field $\sigma_{ij}^{h}$, and the fluctuation field in the stress-free medium in the vicinity of the cavity with imposed tractions $ -\sigma_{ij}^{h} \cdot n_j$  (Fig.3(A)). The second problem is solved with BEM. The formulation remains the same as in the finite-domain case, however now the integrals need to be computed over an infinite periodic surface. Fortunately, FMM can easily be extended to compute this type of integrals: as the multipole expansions are also periodic, it is sufficient to make the FMM tree periodic at all levels \cite{pvfmm}.  This leads us to a solution of the problem depicted in Fig.1(B). Note that within this formulation the boundaries of a periodic cell are not included, as this is not a part of the infinite periodic boundary.

\begin{figure}
	\includegraphics[width=16cm]{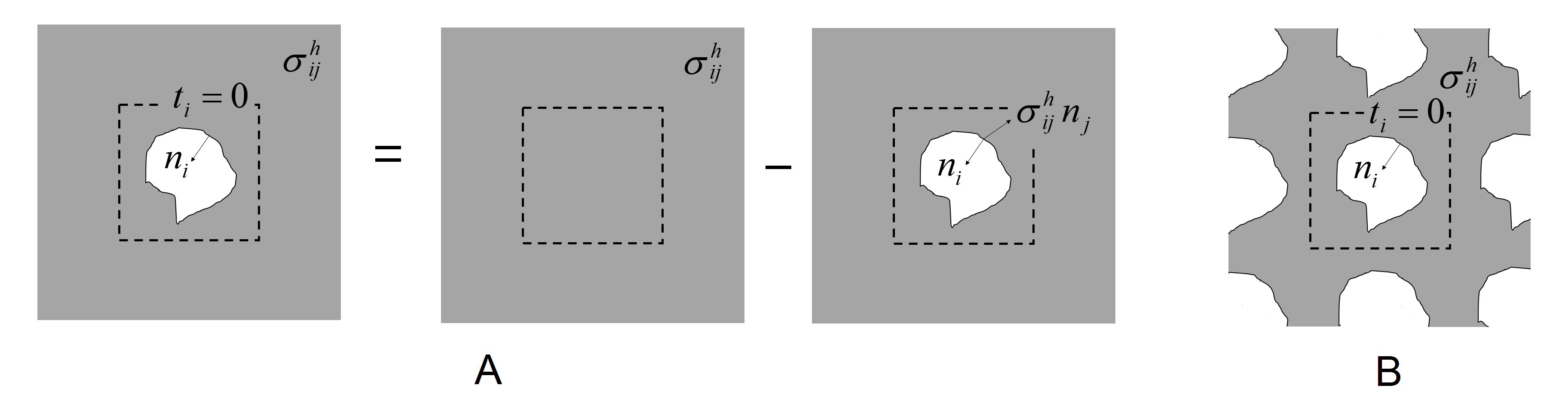}
	\protect\caption{(A) The problem of the traction-free cavity in the medium with otherwise uniform state of stress. (B) Extension to a periodic system of cavities.}
\end{figure}

Using this BVP formulation, we use the standard TD (\ref{deriv}) with the imposed hydrostatic stress, in order to obtain the periodic cell with the maximum bulk modulus. Since we do not impose isotropy constraints, the obtained periodic cell exhibits cubic orthotropic elastic behavior defined by the cubic symmetry of the cell. Clearly, our formulation requires the initial guess on the shape, since we cannot start the optimization process with the homogeneous field of topological derivatives and absent boundaries. As such, we use the configuration with the spherical hole placed in the center of the cell. The obtained optimal configuration may depend on the initial guess (in our particular case, on the radius of the cavity $R_c$).

\begin{figure} 
	\includegraphics[width=16cm]{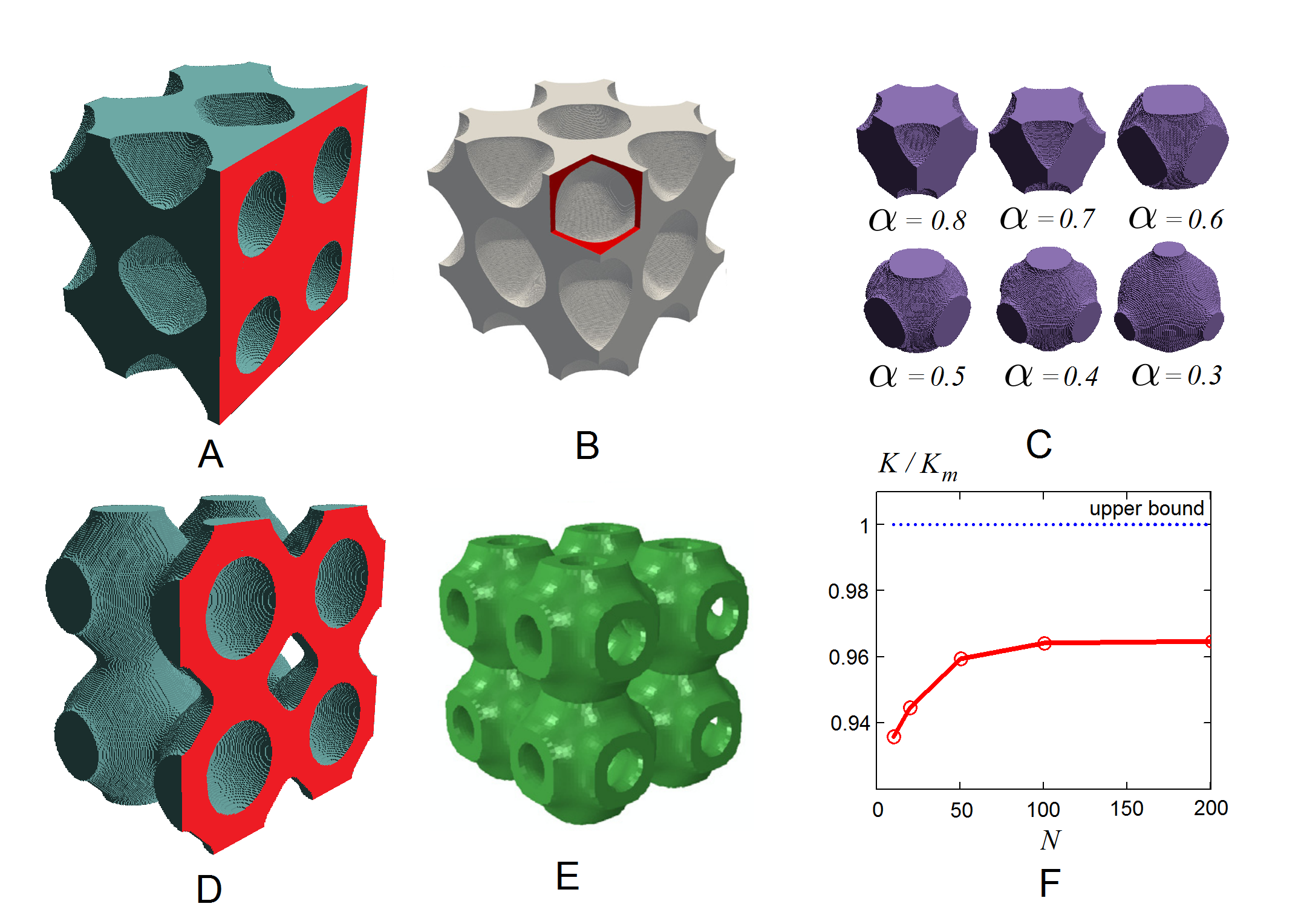}
	\protect\caption{ Periodic cells of an elastic metamaterial with maximum bulk modulus ($2 \times 2 \times 2$ cells are shown). (A) Single-iteration BEM solution obtained for $\alpha_c = 0.7$ ($100^3$ voxels, $R_c = 0.2$). (B) Maximum bulk modulus cell with the material volume fraction $0.7$ obtained in \cite{petsc_fastest} (C) Single-iteration BEM solution obtained for $\alpha_c = 0.4$ ($100^3$ voxels, $R_c = 0.3$). Maximum bulk modulus cell with volume fraction $0.4$ obtained in \cite{Huang2011}.}
	
\label{bulk}
\end{figure}

Figures \ref{bulk} (A), (D) show the results of our optimization. The models depicted in \ref{bulk} (A), (D) consisted of $100^3$ voxels,  had $3 \cdot 10^6$ volume DOFs  and $6.2 \cdot 10^6$ surface DOFs at the final iteration of the optimization process. Unlike the cases of mixed BVPs and concentrated force loadings, for Neumann formulation of a periodic cell problem GMRES convergence process took only $5-10$ iterations, which resulted to an excellent computation time, less than a minute per iteration on 128 cores. Figures \ref{bulk} (A), (D) depict the results of a single-iteration convergence. As we can see, these configurations are in good aggreement with the solutions of the same problem obtained with a lot more sophisticated and computationally  expensive FEM-fased approaches - method of moving asymptotes (MMA) (Fig. \ref{bulk} (B)) and bidirectional evolutionary structure optimization (BESO) (Fig. \ref{bulk} (E)). Furthermore, the obtained configurations' bulk moduli are within 95\% of Hashin-Shtrickmann upper bound \cite{Hashin} for a material-void composite.

Figure \ref{bulk} (C) demonstrates the evolution of the single-iteration optimal solutions for different volume fractions($R_c = 0.3$). Figure \ref{bulk} (E) gives the dependence of the resulting bulk modulus of the cell as a function of the number of voxels along the cube of the model. As we can see, there is clear convergence of the energy cost functional (and the corresponding bulk modulus ) with the refinement of the mesh.  
  
As we could see, a simplistic single-iteration approach for topology optimization appeared surprisingly efficient. It is of particular interest for the problems of topological optimization of periodic cells of metamaterials, since in this particular class of problems 1) one is interested in possibly simple shapes and topologies of the microstructure, additional levels of microstrucutre are undesirable; 2) one often need to perform wide parametric studies of the shapes and topologies that depend on the initial guess, and it is crucial to have a fast and scalable solver, that can handle thousands of optimization cycles.

\begin{figure}  \label{3Dprint}
	\includegraphics[width=16cm]{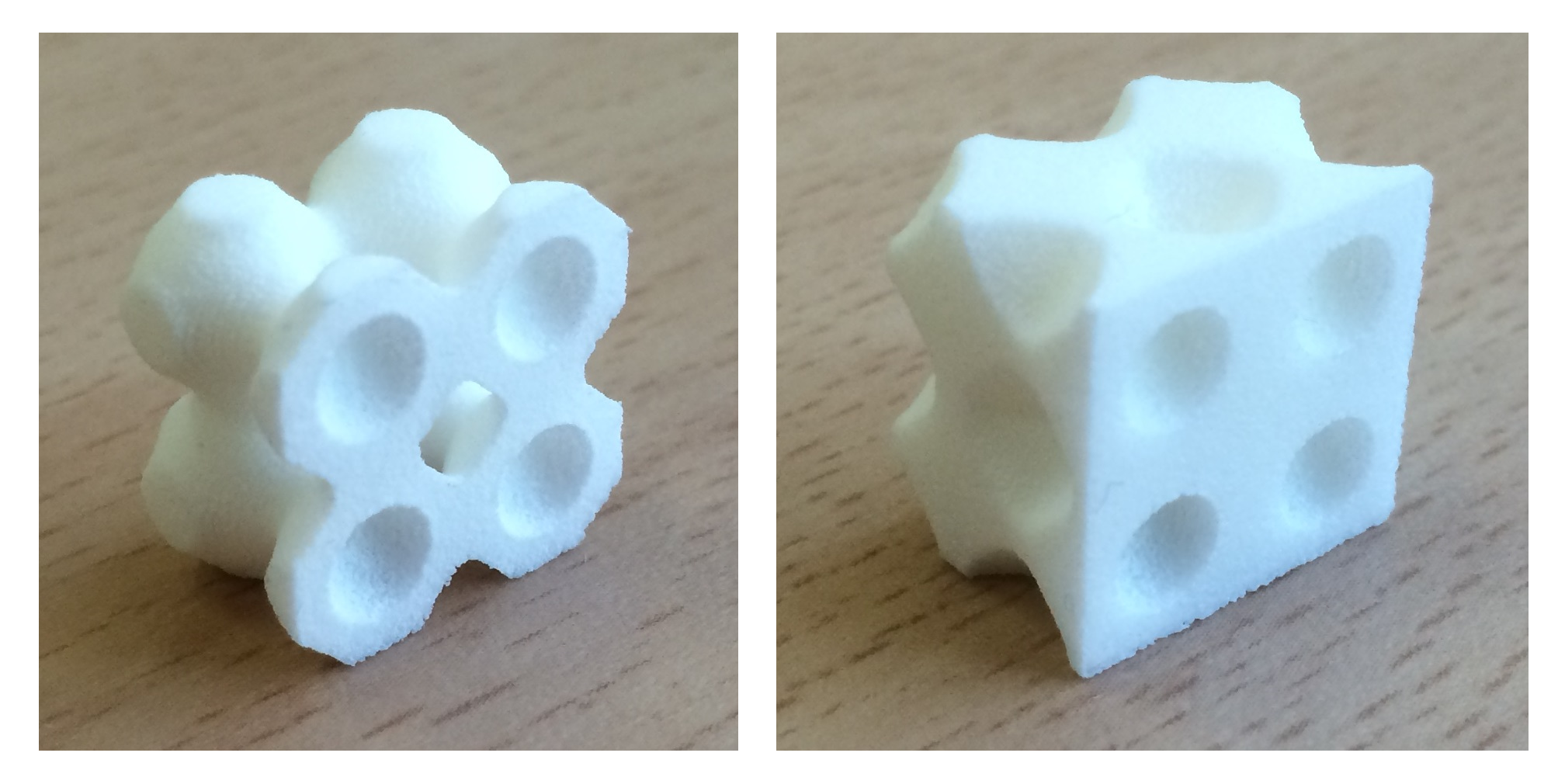}
	\protect\caption{ 3D printed prototypes of the maximum bulk modulus microstructures}
\end{figure}

It worth noting here that the optimal configurations can be easily saved in STL (table of triangle's vertices and normals) format, that can be immediately used for additive manufacturing of optimal structures. Figure \ref{3Dprint} presents the maximum bulk modulus microstructures shown in Fig. \ref{bulk} (A,D), rendered in polyamide plastic using selective laser sintering technology.

\section{Discussion and future work} \label{disc}

Our paper presents the first  scalable implementation of a three-dimensional BEM-based topology optimization algorithm. It is therefore interesting to compare it with the state of the art FEM implementations in terms of robustness, performance, and specific features. 

One of the unwanted yet ubiquitous features of FEM topology optimization techniques are checkerboard patterns \cite{checkerboard}. Checkerboards appear because  of high stiffness of the checkerboard pattern in finite-element discretization, in  comparison with a continuous density distribution with the same total mass.  Consistently with previous  work  \cite{BEM2DMarczak, BEM3D} we  do not observe checkerboard patterns in our formulation. In all our simulations, in spite the low-order approximation and absence of explicit regularization, we did not observe anything similar to typical FEM  checkerboards.

Another important feature of FEM optimization techniques is their inherent dependence on the volume mesh. In absence of regularization the cost functional achieved in the FEM topology optimization process is often heavily dependent on the level of grid refinement \cite{Meshdep}. Our simulations do not show  much dependence of the structure obtained in optimization process on the level of grid refinement (Fig. \ref{bulk}, (F)). As has been mentioned above, the simulation result strongly depends on the level of threshold parameter $C$, which regularizes the problem and defines the structure and the corresponding value of the cost functional. 

Since the "hard-kill" greedy algorithm of material removal is used, we can not guarantee that the solutions found within our approach are indeed globally optimal. However, both qualitative shapes and the functional values demonstrate that our optimization algorithm finds the solutions that are close to what is found with FEM homogenisation techniques.

Our code provides the single node performance and parallelization comparable to the best available FEM codes. For example, the problem of finding the optimal periodic cell for the highest bulk modulus was addressed in \cite{petsc_fastest}. The solution obtained in \cite{petsc_fastest} is very similar to ours. The discretization used $288^3$ design degrees of freedom, and the solution took $60$ second per iteration on $240$ cores. However, within the method of moving asymptotes used in  \cite{petsc_fastest} hundreds of iterations are required to achieve the final extremal structure. Our single-iteration solution for $200^3$ is obtained in just two minutes on $128$ cores. Based on these results we can claim that our optimization technique is at least of comparable performance with state of the art FEM techniques. This performance, as well as the quality and the robustness of solutions, can be further improved. In the remainder of the section we discuss possible directions of further development of our technique.

The first and straightforward improvement is the introduction of higher order boundary elements. Piecewise-constant elements were chosen in our work because of simplicity of parallelization - such approximation scheme does not require interprocess communications beyond those already implemented in PVFMM. However, for better convergence and quality of optimal solutions it is necessary to improve the order of approximation.

	\begin{figure} \label{future}
		\includegraphics[width=16cm]{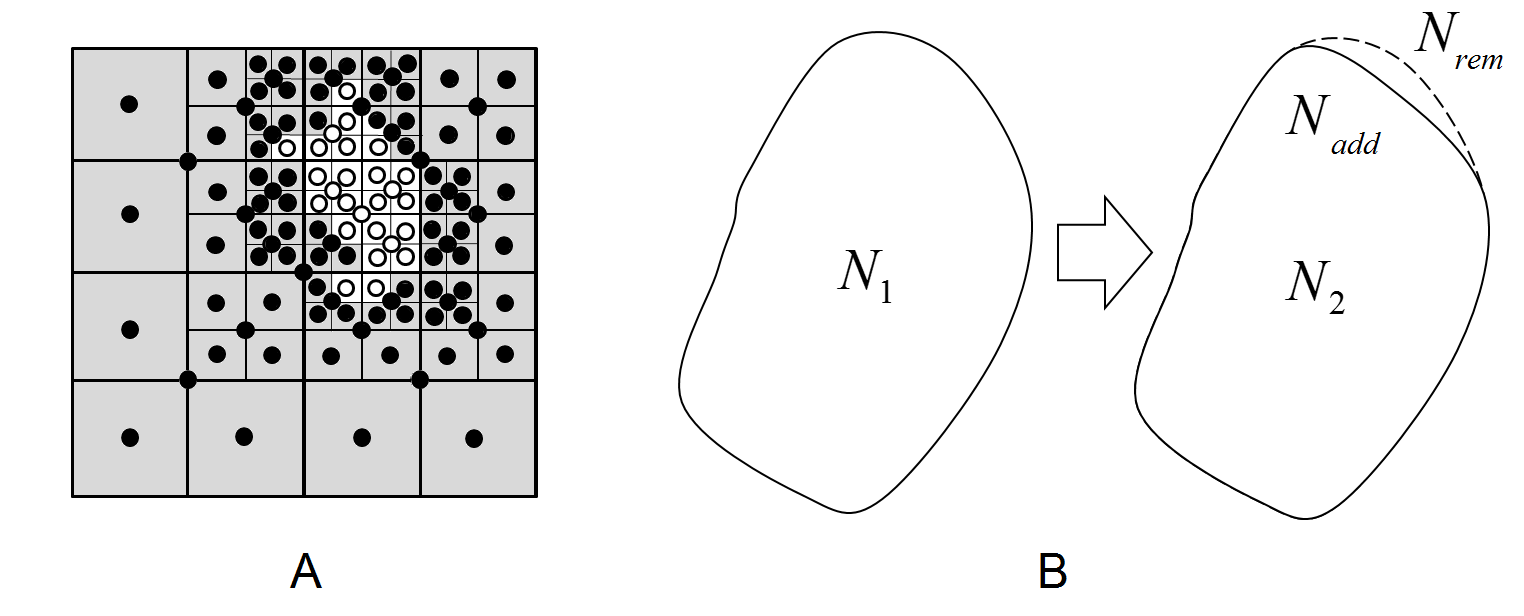}
		\protect\caption{(A) Adaptive sampling of the topological derivatives (B) Fast update of the boundary solution and the fields inside using the solution from the previous iteration}
	\end{figure}

 The second important direction is improving the surface mesh generation. In this work we use simple approach to mesh generation by voxel sides with subsequent Laplace smoothing of the resulting configuration. Such mesh generator provides sufficient quality of the mesh for the first-order convergent piecewise-constant elements. However, in case of higher-order elements a better surface approximation is required. The marching cubes algorithm and its generalizations \cite{mcubes,cazhdan,dual} are more appropriate in these cases.

 In this work we demonstrated the application of our code to the problems dealing with linear elasticity fundamental solutions. However, due to kernel independence our approach can be straightforwardly generalized onto a number of other non-oscillatoric kernels, including, for example, Laplace and  Stokes kernels. 

Presented method does not fully address the issue of computation precision loss for the fields inside the domain, when the point of interest approaches the boundary. For our surface discretization, this distance is never below one half of the voxel size, and thus this problem can be handled by boundary refinement.  However, in a more general case, when marching cubes-type algorithm is used, the target evaluation points may be arbitrarily close to the surface, and an interpolation method is needed (e.g., \cite{ying2006high}).

  In the examples presented above we have used the uniform volume grid. However, domain computations can be performed on an adaptive grid, reducing therefore the computational complexity to linearly proportional with respect to number of surface discretization elements. Two-scale adaptive grid in this case is adjusted to coarse and fine lengthscales: the coarse one is defined by the minimal size of the topological feature we would like to detect, and fine scale is defined by the shape features we whould like to resolve for this topology. The procedure of calculation of the topological derivatives is therefore performed stepwise, from coarse to fine level (Fig. (\ref{future}) (A)). The simplest criteria for grid refinement - different values of threshold function for two neighboring cells. For such an adaptive computation scheme domain field computation would require $O(d \cdot N_s)$ operations, where $d$ is the number of levels of the domain grid refinement. Such an adaptive scheme is therefore useful only if the number of volume points $N_v$ is significantly larger than the product  $d \cdot N_s$, which is the case only for the models that are much larger than those considered in the current work.

As it was noted in the examples section, GMRES convergence process appears to be slow for the case of the mixed boundary condition problems and problems in Neumann formulations with concentrated forces. In the future this shortcoming should be addressed with the better choice of integral formulation or the appropriate preconditioner. 

As we could see in the considered examples, for a wide class of problems the desirable topology is found in a single iteration, whereas subsequent iterations of the topology optimization work just as a brute-force shape optimization. This leads us to a conclusion that the efficient BEM-based algorithm of the topological-shape optimization should combine the single iteration of the topology optimization followed by shape optimization with a lot more time- and memory-efficient shape optimization formalisms based on shape derivatives and shape gradients \cite{NovotnySokolovskyBook}. 

In our earlier works \cite{Ostanin,Ostanin_wit} we have demonstrated that if the change in boundary configuration at every iteration is relatively small (Fig. (\ref{future}) (B)), one can use fast update techniques for the volume and surface solutions, which would be faster than full re-computation of the BVP. The suggested technique of the fast update of the surface solution is based on Shur complements\cite{Inverse}, whereas the technique for the fast update of the field of the topological derivatives was based on the superposition of the scalar products of the partially known influence coefficients and old/new boundary solution. These techniques were described in a context of small models, and full system matrix representation. The development of the analogous tools in a context of FMM factorization would become a significant advance of BEM-based topology optimization.

\section{Related work}

In this section we provide a brief literature survey that puts our work in a context of recent achievements in the field. A number of efficient optimization techniques has been developed during the last few decades. They can be divided into two broad sets - various composite material homogenization techniques \cite{Bendsoe1988, AllaireHomog,Topopt}, that optimize the distribution of the material density, which is then thresholded, and binary optimization techniques, that prohibit intermediate material density at the optimization stage\cite{Novotny3Delast,Sokol3D,NovotnySokolovskyBook}. The virtue of the first kind of approaches - wider search space, that in many cases facilitates finding better designs and problem convexification \cite{Topopt}. The strength of the approaches of the second kind is a complete description of the optimization problem given in terms of the surface of the domain, so the estimates in the functional value and gradients computed at each step of the optimization process correspond to the actual material distribution, not a smoothed version of it used in approaches of the first type. An additional benefit that we exploited in this paper is that for linear elasticity (or any other linear PDE), the solution can be obtained using a boundary integral formulation.
Until recently, BEM techniques were not used for topological optimization problems. These were, however, applied to solve inverse scattering problems in elastodynamics \cite{Bonnet}, which are related to topology/shape optimization.
First applications of BEM to topological optimization of elastic structures were presented in \cite{BEM2DMarczak,BEM3D}. These early works demonstrated conceptual applicability of BEM in combination with a hard-kill algorithm of material removal to the problems of topology optimization. The first applications of algebraically accelerated BEM to two- and tree- dimensional problems of elasticity were presented in our papers \cite{Ostanin,Ostanin_wit}. In these works we used Shur complements \cite{Inverse} for fast updates of the BVP solution \cite{Ostanin}, and $\mathcal{H}^2$-matrices for fast solutions of the BVP \cite{Ostanin_wit}. Nonetheless, these implementations were neither truly scalable nor parallelizable.
In this work we have presented the first scalable parallel realization of the BEM-based topology optimization algorithm, suitable for topological optimization with millions of degrees of freedom.

\section{Conclusion}

In this work we presented the first technique for large-scale topological-shape optimization based on FMM-accelerated BEM. The approach uses direct boundary element formulation and the kernel-independent fast multipoole method. The method utilizes voxel representation of the domain and iterative isosurface extraction of topological derivatives. The obtained approach is free of the typical shortcomings of FEM-based techniques, such as checkerbord instabilities and mesh dependent optimization results. The efficiency of the proposed technique was illustrated on the examples of minimum compliance structural optimization, as well as the optimization of the periodic cell of the material for the desirable tensor of elasticity. 

\section{Acknowledgements} 

Authors express their deep gratitude to Dhairya Malhotra, for his helpful comments and assistance. 
Authors gratefully acknowledge the financial support from Russian National Foundation under the grant №15-11-00033. I.O acknowledges the financial support from the Russian Foundation of Basic Research under grant 16-31-60100.


\bibliographystyle{unsrtnat}
\bibliography{manuscript}

\end{document}